\documentclass{amsart}

\usepackage{amssymb}

\usepackage{graphicx}

\newtheorem{theorem}{Theorem}[section]
\newtheorem{lemma}[theorem]{Lemma}
\newtheorem{proposition}[theorem]{Proposition}
\newtheorem{corollary}[theorem]{Corollary}
\theoremstyle{definition}

\theoremstyle{remark}

\numberwithin{equation}{section}

\begin{document}

\title[The Loewner driving function of trajectory arcs]{The Loewner driving function of trajectory arcs of quadratic differentials}


\author[J. Tsai]{Jonathan Tsai}
\address{Department of Mathematics \\ Chinese University of Hong Kong \\ Shatin \\ New Territories \\ Hong Kong}
\email{jtsai@math.cuhk.edu.hk}
\thanks{}


\begin{abstract}
We obtain a first order differential equation for the driving function of the chordal Loewner differential equation in the case where the domain is slit by a curve which is a trajectory arc of certain quadratic differentials. In particular this includes the case when the curve is a path on the square, triangle or hexagonal lattice in the upper half-plane or, indeed, in any domain with boundary on the lattice. We also demonstrate how we use this to calculate the driving function numerically. Equivalent results for other variants of the Loewner differential equation are also obtained: Multiple slits in the chordal Loewner differential equation and the radial Loewner differential equation. The proof of our theorem uses a generalization of Schwarz-Christoffel mapping to domains bounded by trajectory arcs of rotations of a given quadratic differential.
\end{abstract}

\maketitle


\bibliographystyle{amsplain}
\section*{Introduction}
Suppose that $\mathbb{H}=\{z\in\mathbb{C}:\mathrm{Im}(z)>0\}$ is the upper half-plane and $\gamma:[0,T)\mapsto\overline{\mathbb{H}}$ is a simple Jordan curve with $\gamma(0)\in\mathbb{R}$ and $\gamma(0,T)=\{\gamma(t):t\in(0,T)\}\subset\mathbb{H}$. Then for each $t\in(0,T)$,
\[H_{t}=\mathbb{H}\setminus \gamma(0,t]\]
is a simply-connected domain and hence by the Riemann mapping
theorem, we can find a conformal map $f_{t}$ of $\mathbb{H}$ onto
$H_{t}$. Moreover, we can require that
$f_{t}$ has series expansion \[f_{t}(z) = z- \frac{C(t)}{z} +
O\left(\frac{1}{z^2}\right) \text{ as } z\rightarrow \infty.\]
Normalized in this way $f_{t}$ is unique and is said to be
\emph{hydrodynamically normalized}. The function $C(t)$ is positive, continuous and
strictly increasing: it is called the \emph{half-plane capacity} of
$\gamma(0,t]$. Thus we can reparameterize $\gamma$ such that $C(t)=2t$ for
all $t$, we will call this \textit{parameterization by half-plane
capacity}. With this normalization and parameterization, the
function $f_{t}$ satisfies the differential equation (where $f_{t}'$ denotes differentiation with respect to $z$ and $\dot{f}_{t}$ denotes differentiation with respect to $t$):
\begin{equation}\dot{f}_{t}(z)=-\frac{2f_{t}'(z)}{z-\xi(t)}, \label{hatton}\end{equation}
where $\xi(t)=f_{t}^{-1}(\gamma(t))$ is a continuous real-valued function. This
is the \emph{chordal Loewner differential equation}; $\xi(t)$ is
called the \emph{driving function} of the slit $\gamma$. The
converse is also true: given a measurable function $\xi$, the
differential equation (\ref{hatton}) with initial condition
$f_{0}(z)\equiv z$ has solution $f_{t}$ which is a conformal map
from $\mathbb{H}$ into itself (although $f(\mathbb{H})$ is not
necessarily a slit domain). Chapter 3 and 4 of \cite{Law05} gives full
details of this construction.

\begin{figure}
\scalebox{0.4}{\includegraphics{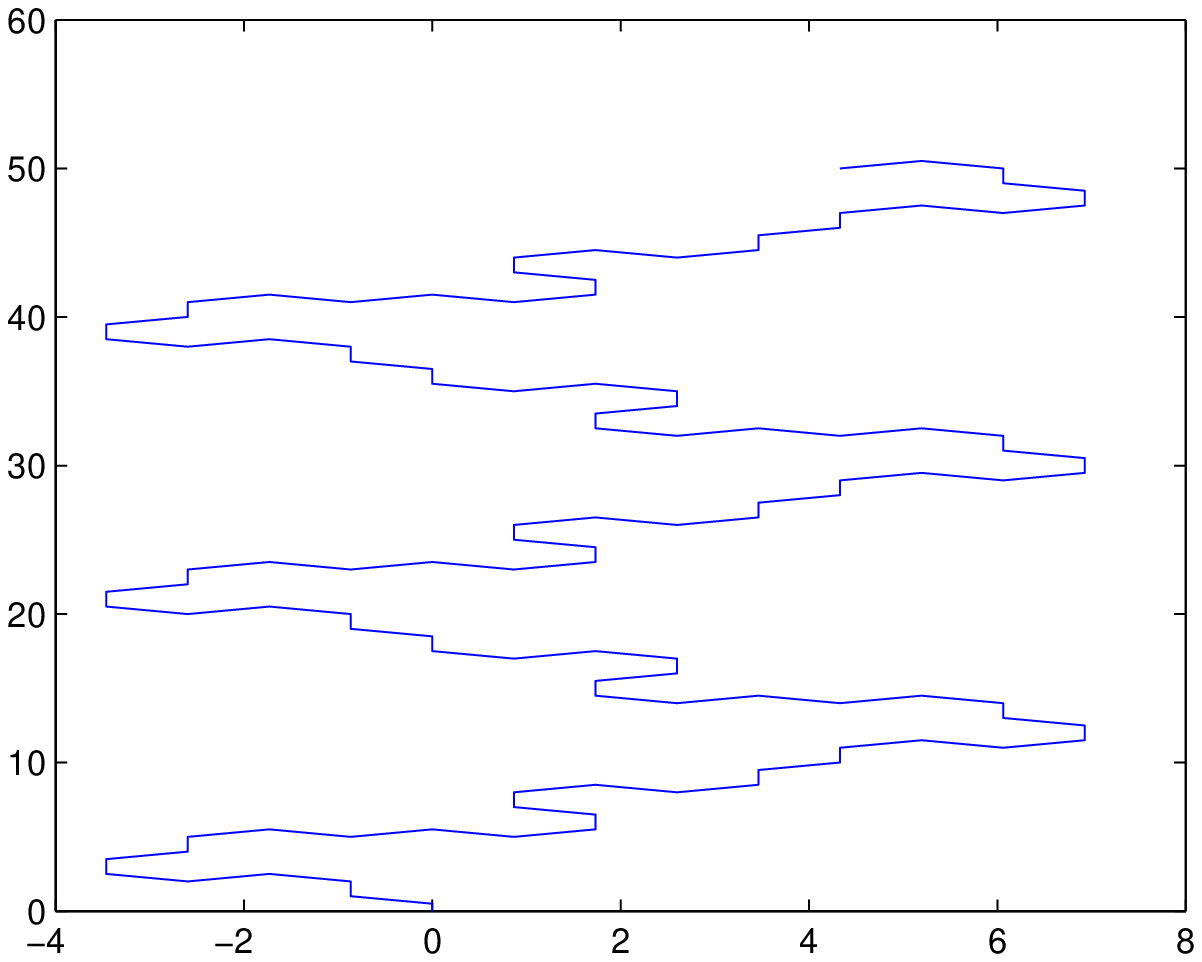}\includegraphics{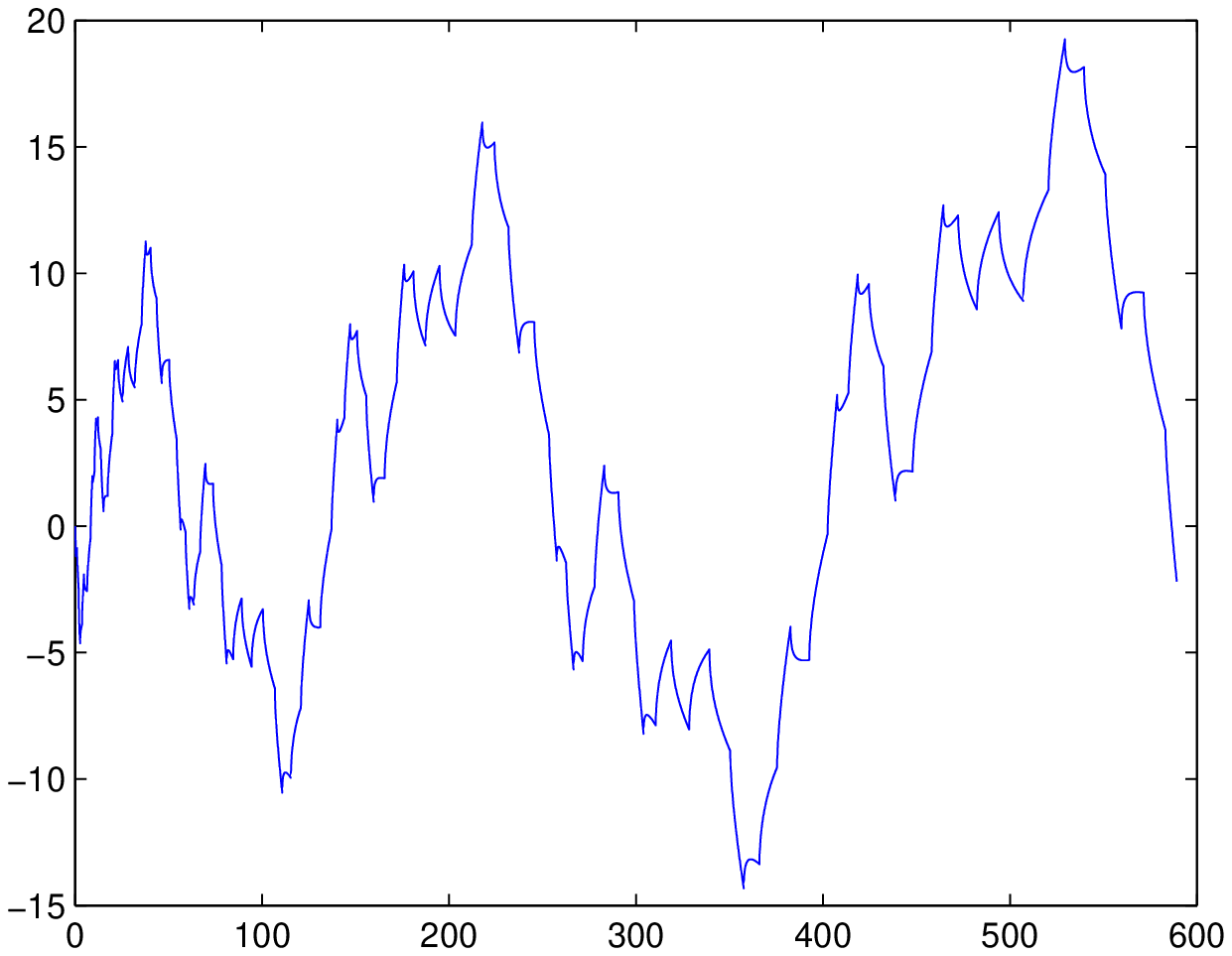}}
\caption{A path on the hexagonal lattice on the upper half-plane
(left) and a plot of its driving function on the $y$-axis against
time on the $x$-axis (right).} \label{fig1}
\scalebox{0.4}{\includegraphics{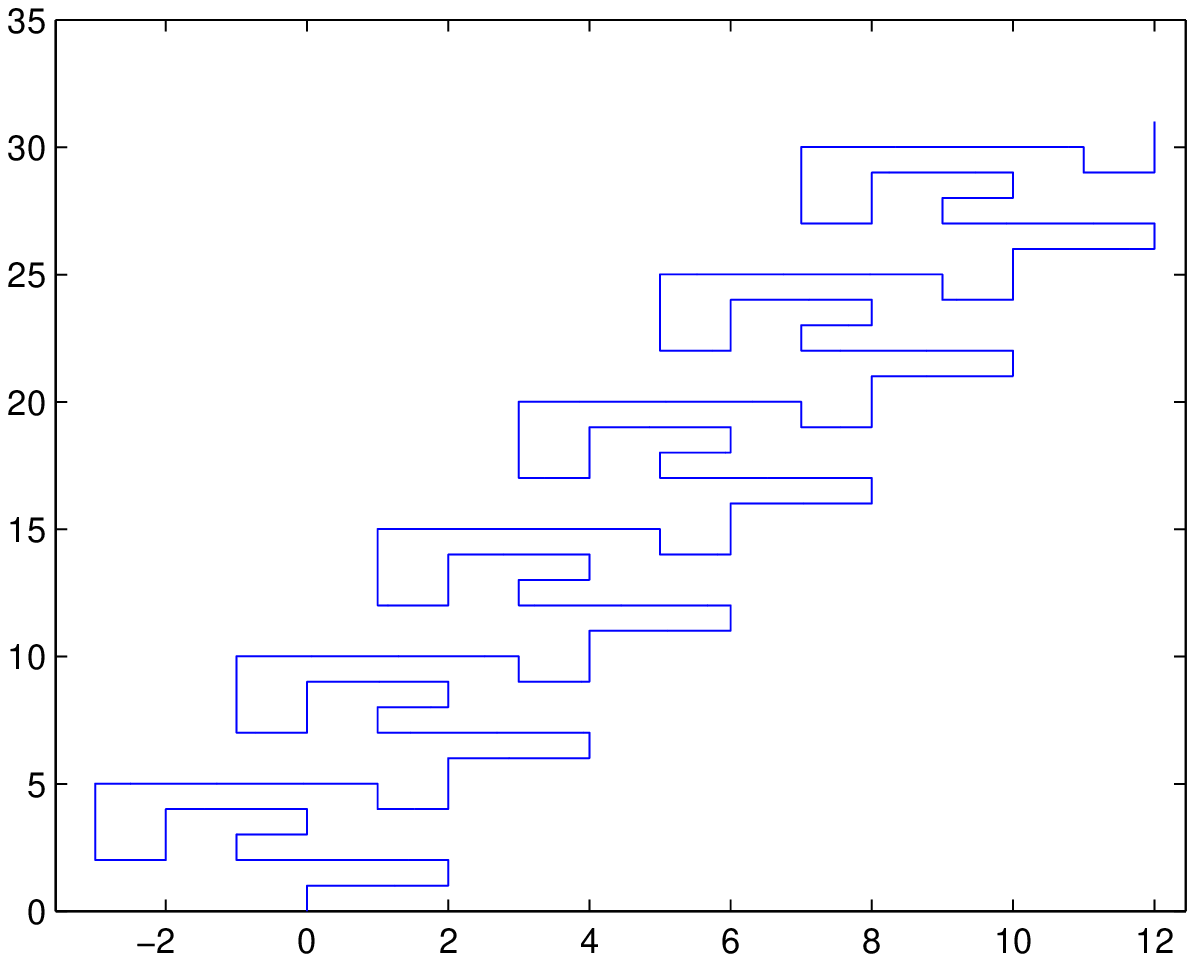}\includegraphics{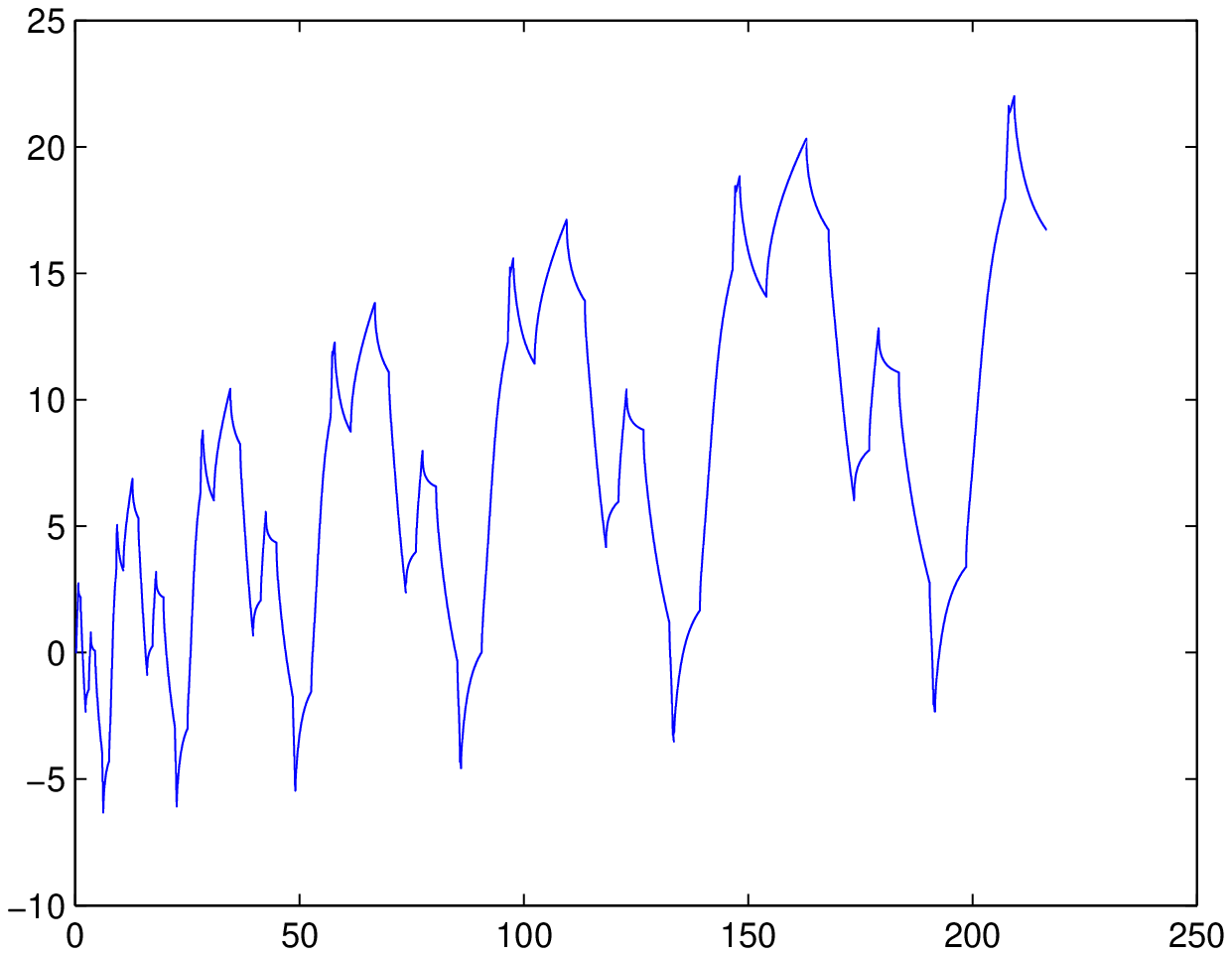}}
\caption{A path on the square lattice on the upper half-plane (left)
and a plot of its driving function on the $y$-axis against time on
the $x$-axis (right).} \label{fig3}
\end{figure}

Since Schramm's discovery of stochastic Loewner evolution in 1999
(see \cite{Schramm00}), there has been huge interest in the chordal
Loewner differential equation and its variants. But the relationship
between  the slit in $\mathbb{H}$ and its resulting driving function
is not well understood. There are a few papers that relate the
behaviour of the slit with the behaviour of the driving function
e.g. \cite{MR05},\cite{Lind05}; also, the paper \cite{KNK04}
calculates the slit arising from a few driving functions. In this
paper, we will obtain a first order differential equation for $\xi$
(which we can then solve numerically) that
allows us to calculate the driving function $\xi$ in the case where
the curve $\gamma$ is a trajectory arc of a certain type of quadratic
differential. We will show that this includes, for example, the
case when $\gamma$ is a path on the square/triangle/hexagonal
lattice in the upper half-plane or indeed, in any domain whose boundary
lies on such a lattice. So for example, Figure \ref{fig1} plots the driving
function of a path on the hexagonal lattice in the upper half-plane and Figure \ref{fig3}
plots the driving function of a path on the square lattice in the upper half-plane.

We also note that we can obtain equivalent results for other variants of the Loewner differential equation for example, in the radial version or with multiple slits. We will discuss this in the paper as well.

The proof of our formulae uses a generalization of Schwarz-Christoffel mapping to domains bounded by trajectory arcs of rotations of a given quadratic differential.

We also mention that, currently, the common method used to find the
driving function of a given slit is to use the Zipper algorithm
discovered independently by D. E. Marshall and R. K\"{u}hnau to
approximate the function $f_{t}$ which can then be used to determine the
driving function. The Zipper algorithm can be viewed as a discrete
version of the Loewner differential equation and hence is well
suited to studying growth processes. It also has the advantage of
being very fast. See \cite{MR06} and \cite{Ken07}.

\section{Main results}
To state our main results, we have to provide some background in the
theory of quadratic differentials. Note that not all the terms used
here are standard in the literature. See Chapter 8 of \cite{Pom75}
and \cite{Streb80} for more details. A \textit{quadratic
differential} on a domain $D\subset\widehat{\mathbb{C}}=\mathbb{C}\cup\{\infty\}$ is the formal expression
\[Q(z)dz^2,\]
where $Q(z)$ is a meromorphic function on $D$. Then for $\omega\in D$ with $\omega\neq \infty$, $Q(z)$ has Laurent series expansion about $\omega$,
\[Q(z)=\sum_{k=n}^{\infty} a_{k}(z-\omega)^{k}\]
for some $n>-\infty$ with $a_{n}\neq 0$. Then we define the \emph{degree of $\omega$ with
respect to $Q(z)dz^2$}, $\deg_{Q}(\omega)$, to be equal to $n$.

If $\infty\in D$, then  near $\infty$, $Q$ has Laurent series expansion given by
\[Q(z)=\sum_{k=m}^{\infty}b_{k}z^{-k},\]
then we define the \emph{degree of $\infty$ with respect to $Q(z)dz^2$}, $\deg_{Q}(\infty)$ to be equal to $m-4$. The ``4'' in the definition ensures that the degree is conformally invariant in a way which we will make precise later. Then $\omega\in D$ is:
\begin{itemize}
 \item a \textit{zero} of $Q(z)dz^2$ if $\deg_{Q}(\omega)>0$.
\item a \textit{pole} of $Q(z)dz^2$ if $\deg_{Q}(\omega)<0$.
\item an \textit{ordinary point} of $Q(z)dz^2$ if $\deg_{Q}(\omega)=0$.
\end{itemize}

A \emph{trajectory
arc} of $Q(z)dz^2$ is a curve $\gamma:(a,b)\mapsto D$ that does not meet any zeroes and poles of $Q(z)dz^2$ and satisfies
\[Q(\gamma(t))\dot{\gamma}(t)^2>0 \text{ for all } t\in (a,b).\]
For $\theta\in[0,\pi)$, a \emph{$\theta$-trajectory arc} of $Q(z)dz^2$ is a curve $\gamma:(a,b)\mapsto D$ that satisfies
\[\arg[Q(\gamma(t))\dot{\gamma}(t)^2] =2\theta  \text{ for all } t\in (a,b).\]
Then $\gamma$ is a $\theta$-trajectory arc of $Q(z)dz^2$ if and only if it is a trajectory arc of $e^{-2i\theta}Q(z)dz^2$. Hence, a $0$-trajectory arc is simply a trajectory arc and we call a $\pi/2$-trajectory arc an \emph{orthogonal trajectory arc}. It is clear that these definitions are invariant under reparameterization of $\gamma$ so we will often call the point set of $\gamma$ a trajectory arc or $\theta$-trajectory arc. We call a maximal trajectory arc a \emph{trajectory} and similarly, a maximal $\theta$-trajectory arc is called a \emph{$\theta$-trajectory}.
For example, if we consider the quadratic differential $1dz^2$ in $\mathbb{C}$, then the $\theta$-trajectories are the straight lines with gradient $\exp(2\theta)$.

We now consider a special type of quadratic differential: Let $D$ be a domain with piecewise analytic boundary. A \emph{K\"{u}hnau quadratic differential} is a quadratic differential, $Q(z)dz^2$, on  $D$ satisfying the following two properties:
\subsection*{Definition:} Let $D$ be a domain with piecewise analytic boundary.  A \emph{K\"{u}hnau quadratic differential} is a quadratic differential, $Q(z)dz^2$, on  $D$ satisfying the following two properties:
\begin{enumerate}
\item We can write
\[\partial D=\bigcup_{j=1}^{n}\overline{\Gamma}_{j}\]
such that each $\Gamma_{j}$ is an open analytic arc with $\Gamma_{k}\cap\Gamma_{j}=\emptyset$ for $k\neq j$ and moreover, $Q(z)$ extends continuously to each $\Gamma_{j}$ and $\mathrm{arg}[Q(z)dz^2]$ is constant on each $\Gamma_{j}$ i.e. each $\Gamma_{j}$ is a $\theta_{j}$-trajectory arc for some $\theta_{j}\in[0,\pi)$.
\item At $z\in \overline{\Gamma}_{k}\cap \overline{\Gamma}_{j}$ for all $j\neq k$, there are either only finitely many direction from which trajectories approach the point $z$ or if there are infinitely many directions from which trajectories approach the point $z$, then for each such direction, there is only one trajectory that approaches $z$ at this direction.
\end{enumerate}
These quadratic differentials are studied by K\"{u}hnau in \cite{Kuh67} where he applies them to the study of certain Gr\"{o}tzsch-style extremal problems

Property (i) above, also implies that $D$ is locally connected. Thus each prime end of $D$ corresponds to a unique point in $\partial D$ (see \cite[p. 27]{Pom92}). If, in addition, a point on $\partial D$ corresponds to a unique prime end, then we make no distinction between the two. Let $z$ be a prime end of $D$.  Then we have 2 cases: Either $z\in\Gamma_{j}$ for some $j=1,\ldots,n$; or there exist exactly 2 of the $(\Gamma_{j})$, that end at the prime end $z$. In the latter case, we will denote $z$ by $z_{k}$ and assume that $\Gamma_{k}$, a $\theta_{k}$-trajectory arc, and $\Gamma_{k-1}$, a $\theta_{k-1}$-trajectory arc, are the only 2 arcs that end at $z_{k}$.
Then we can define \textit{the degree of $z_{k}$ in $D$ with respect to $Q(z)dz^2$}, $\mathrm{deg}_{D,Q}(z_{k})$, as follows:
\[\mathrm{deg}_{D,Q}(z_{k})=\left\{\begin{array}{ll}2[|\theta_{k}-\theta_{k-1}|/\pi +J_{k}-1] & \text{ if }\theta_{k}\neq\theta_{k-1},\\
2J_{k}&\text{ if }
\theta_{k}=\theta_{k-1},
\end{array}\right.\] where $J_{k}$ is
the number of trajectories of $Q(z)dz^2$ inside $D$ that end at the
prime end $z_{k}$. If $J_{k}$ is infinite, then the degree is not defined. Then for prime ends $z$ such that $z\in\Gamma_{j}$ for some $j=1,\ldots,n$, we define
\[\mathrm{deg}_{D,Q}(z)=0.\]
Although the motivation for this definition currently seems unclear, we will see that this indeed generalizes the concept of degree to points on the boundary. In particular, we will show that for $x\in\partial \mathbb{H}$, if $\mathrm{deg}_{\mathbb{H},Q}(x)\in\mathbb{Z}$, then $Q$ can be extended to a meromorphic function in a neighbourhood of $x$ with
\[\mathrm{deg}_{\mathbb{H},Q}(x)=\mathrm{deg}_{Q}(x).\]

We then have the following theorem on K\"{u}hnau quadratic differentials in $\mathbb{H}$:
\begin{theorem}\label{QDSC}
Suppose that $Q(z)dz^2$ is a K\"{u}hnau quadratic differential on
$\mathbb{H}$. Then we have
\[Q(z)=R\left(\prod_{j=1}^{n}(z-\zeta_{j})^{\lambda_{j}}\right)\]
for some constant $R\neq 0$, $\zeta_{j}\in \mathbb{C}$, $\lambda_{j}\in
\mathbb{R}$ for $j=1,\ldots,n$.\end{theorem}
This theorem can be viewed as a generalization of the
Schwarz-Christoffel formula to domains bounded by
$\theta_{k}$-trajectory arcs of a given quadratic differential.

\ \\
We then have the following theorem on the Loewner driving function of a
$\phi$-trajectory arc of a K\"{u}hnau quadratic differential
$Q(z)dz^2$ that starts at a point $\xi_{0}\in\mathbb{R}$ with
$\deg_{\mathbb{H},Q}(\xi_{0})=N\in \{0,1,\ldots\}$.
\begin{theorem}\label{QDdriving}
Suppose that $Q(z)dz^2$ is a K\"{u}hnau quadratic differential on
$\mathbb{H}$ such that there is a point $\xi_{0}\in\mathbb{R}$ with
$\deg_{\mathbb{H},Q}(\xi_{0})=N\in \{0,1,\ldots\}$; then we
have
\[Q(w)=(w-\xi_{0})^N\left(\prod_{j=1}^{n}(w-a_{j})^{\alpha_{j}}   \right),\]
where $a_{j}\in\mathbb{C}$ and $\alpha_{j}\in\mathbb{R}$.
Let $\gamma:[0,T)\mapsto \partial{\mathbb{H}}$ be a simple curve such that $\gamma(0)=\xi_{0}$, $\gamma(0,T)\subset\mathbb{H}$ and $\gamma(0,T)$ is a $\phi$-trajectory arc of
$Q(z)dz^{2}$ ($\phi\in[0,\pi)$) that is parameterized by half-plane
capacity. Suppose that the functions $f_{t}$ maps $\mathbb{H}$ conformally onto $\mathbb{H}\setminus\gamma(0,t]$ and are hydrodynamically normalized. Then for $t\in(0,T)$
\begin{equation}2\xi(t) = -\mu^{-}C^{-}(t) -\mu^{+}C^{+}(t) - \left(\sum_{j=1}^{n} \alpha_{j}A_{j}(t)\right) +\Sigma_{0},\label{toney}\end{equation}
and
\begin{equation}\dot{\xi}(t) = -\frac{\mu^{-}}{C^{-}(t)-\xi(t)} -\frac{\mu^{+}}{C^{+}(t)-\xi(t)} - \left(\sum_{j=1}^{n} \frac{\alpha_{j}}{A_{j}(t)-\xi(t)}\right)\label{wright}\end{equation}
with initial condition $\xi(0)=\xi_{0}$. Where the functions $A_{j}(t)$ are defined by
\[A_{j}(t)= f_{t}^{-1}(a_{j}) \text{ for } j=1,\ldots,n,\]
and $C^{+}(t)>C^{-}(t)$ are the two preimages of $\xi_{0}$ under $f_{t}$;
\[\mu^{\pm}=\mathrm{deg}_{\mathbb{H}\setminus \gamma(0,t],Q}(f_{t}(C^{\pm}(t))),\]
 and
\[\Sigma_{0}=\left[N\xi_{0}+  \left(\sum_{k=1}^{n} \alpha_{k}a_{k}\right) \right].\]
\end{theorem}
We can then use Theorem \ref{QDdriving} to find the driving function in the case
when the slit $\gamma$ consists of consecutive
$\theta_{k}$-trajectory arcs of given quadratic differentials. We
will explain how to do this in further detail later. One difficulty
with using Theorem \ref{QDdriving} is that the parameterization is inherently
given in terms of half-plane capacity. This makes it difficult to
calculate the driving function $\xi$ if we do not know anything
about the half-plane capacity of the trajectory arc (which, in
general, is the case). The next theorem will allow us to compare the parametrization with the length of the slit:

\begin{theorem}\label{gammadot}
Suppose that $Q(z)dz^2$, $\gamma$ and $f_{t}$ are as defined in Theorem \ref{QDdriving}. Let
\[\Phi_{t}(z)=\frac{Q(f_{t}(z))f_{t}'(z)^2}{(z-\xi(t))^2}.\]
Then $\gamma$ satisfies
\begin{equation}\dot{\gamma}(t)=-2\sqrt{\frac{\Phi_{t}(\xi(t))}{Q(\gamma(t))}}.\label{erdei}\end{equation}
\end{theorem}

The rest of this paper is organized as follows: In the Section 2, we will state
some basic results from the theory of quadratic differentials and
use them to prove Theorem \ref{QDSC}. Then we will use Theorem \ref{QDSC} to prove Theorems 1.2 and 1.3 in Section 3. In
Section 4 we will discuss how to obtain the driving function
numerically using Theorems 1.2 and 1.3. Finally in Section 5, we
will discuss extensions of Theorem \ref{QDdriving} to the case with multiple
slits as well as to the radial Loewner differential equation.

\section{K\"{u}hnau quadratic differentials and generalized Schwarz-Christoffel mapping}
The aim of this section is to prove Theorem \ref{QDSC}. We will first look at some of the basic results in the theory quadratic differentials that we will need.
\subsection*{Transformation Law}\ \\
Suppose that $f$ is a conformal map from a domain $D_{2}$ onto a domain $D_{1}$ and suppose that $Q_{1}(w)dw^2$ is a quadratic differential on $D_{1}$. If we define
\begin{equation}
Q_{2}(z)\equiv Q_{1}(f(z))f'(z)^2
\label{khan}\end{equation}
then $Q_{2}(z)dz^2$ is a quadratic differential on $D_{2}$. Then, it is clear that $\theta$-trajectory arcs are preserved by this transformation law i.e.
\[\gamma \text{ is a $\theta$-trajectory arc of } Q_{2}(z)dz^2 \Leftrightarrow f\circ\gamma \text{ is a $\theta$-trajectory arc of } Q_{2}(w)dw^2,\]
and also, for $z\in D_{2}$
\[\deg_{Q_{2}}(z)=\deg_{Q_{1}}(f(z)).\]
Hence trajectories and $\deg_{Q}$ are conformally invariant in the above sense.
\ \\

The following lemma tells us that the behaviour of a quadratic differential at a neighbourhood of a point is determined by the degree of that point.

\begin{lemma}[Local behaviour of quadratic differentials]\label{lemloc}
Let $Q(z)dz^2$ be a quadratic differential on a domain $D$. Then for every $\omega\in D$ there is a conformal mapping $w=\phi(z)$ of some neighbourhood of $\omega$ such that
\[Q(z)dz^{2}=\left\{\begin{array}{ll} dw^2 &\text{ if } \deg_{Q}(\omega)=0,  \\ w^n dw^2 &\text{ if } \deg_{Q}(\omega)=n\geq 1, \\w^{-n} dw^2 &\text{ if } \deg_{Q}(\omega)=-n\leq -1 \text{ with $n$ odd}, \\ c^2 w^{-2}dw^2 &\text{ if } \deg_{Q}(\omega)=-2, \\ (w^{-n}+c w^{-1})^2 dw^2 &\text{ if } \deg_{Q}(\omega)=-n\leq -4 \text{ with $n$ even}.\end{array}\right.\]
Here, $c$ is the residue of a branch of $\sqrt{Q(z)}$ at $\omega$.
\end{lemma}
\begin{proof}
See Theorem 8.1 of \cite{Pom75} or Section 6 of \cite{Streb80}.
\end{proof}
So since trajectories are conformally invariant this lemma tells us that the local structure of trajectories around a point $\omega\in D$ is completely determined by $\mathrm{deg}_{Q}(\omega)$ and the converse is true as well.
\begin{lemma}\label{loc}
Suppose that $\omega\in D$ and $\deg_{Q}(\omega)=n$. Then
\begin{enumerate}
\item For $n\geq -1$, there are exactly $n+2$ trajectories of $Q(z)dz^2$ that end at $\omega$ and form equal angles with each other.
\item For $n\leq -3$, there are infinitely many trajectories ending at $\omega$ and moreover, there are $|n|-2$ directions at $\omega$ forming equal angles such that the trajectories approach $\omega$ in these directions.
\item For $n=-2$, the behaviour depends on the value of $c$ (as defined in Lemma \ref{lemloc}).
    \begin{enumerate}
    \item If $c$ is real, then the trajectories are the images of all radial lines under the map $\phi$ defined in Lemma \ref{lemloc}.
    \item If $c$ is purely imaginary, then the trajectories are the images of all concentric circles  under the map $\phi$ defined in Lemma \ref{lemloc}.
    \item If $\mathrm{Re}[c],\mathrm{Im}[c]\neq 0$ then the trajectories are the images of logarithmic spirals under the map $\phi$ defined in Lemma \ref{lemloc}.
    \end{enumerate}
\end{enumerate}
\end{lemma}
\begin{proof} See Section 7 of \cite{Streb80}.
\end{proof}
This lemma shows that it makes sense for us to define $\deg_{D,Q}(x)$, the degree of a point on the boundary, in terms of the trajectories ending at $x$. If $D=\mathbb{H}$ and $Q$ extends to a meromorphic function on a neighbourhood of some $x\in\mathbb{R}\cup\{\infty\}$ with $\deg_{\mathbb{H},Q}(x)$ finite. Then by studying the trajectory structure at $x$, we can see that
 \[\deg_{\mathbb{H},Q}(x)=\deg_{Q}(x).\]
 This is the motivation for defining $\deg_{D,Q}$ in the way we have. The next lemma shows that the $\deg_{D,Q}$ is also conformally invariant:
\begin{lemma}\label{gendeg}
Suppose that $Q(z)dz^2$ is a K\"{u}hnau quadratic differential on a domain $D$ and $f$ is a conformal map of the upper half-plane $\mathbb{H}$ onto $D$. Then the quadratic differential $Q_{f}(w)dw^2$ on $\mathbb{H}$, defined by $Q_{f}(w)=Q(f(w))f'(w)^2$, is also a K\"{u}hnau quadratic differential. Moreover, suppose that
$z\in\partial D$ is a prime end of  $D$. Then
\[\deg_{D,Q}(z)=\deg_{\mathbb{H},Q_{f}}(f^{-1}(z)).\]
\end{lemma}
\begin{proof}
By Carath\'{e}odory's theorem, $f$ extends continuously to $\partial \mathbb{H}$ and by Schwarz's reflection, $f$ extends analytically across $f^{-1}(\Gamma_{k})$ for all $k=1,\ldots,n$. Since $\theta$-trajectory arcs are conformally invariant, this implies that $Q_{f}$ defined by (\ref{khan}) is a K\"{u}hnau quadratic differential on $\mathbb{H}$. Moreover, for all $k=1,\ldots,n$, $f(\Gamma_{k})$ is a $\theta_{k}$-trajectory arc of $Q_{f}(z)dz^2$. Also each point on $\partial\mathbb{H}$ corresponds bijectively to a prime end of $\mathbb{H}$. Hence there is a bijective correspondence between the points of $\partial \mathbb{H}$ and prime ends of $D$. Then
\[\deg_{D,Q}(z)=\deg_{\mathbb{H},Q_{f}}(f^{-1}(z))\]
follows from the conformal invariance of trajectories.
\end{proof}
\subsection*{Reflection across trajectories}\ \\
Suppose that $D$ is a domain such that $\Gamma\subset\partial D$ is
an open interval in $\mathbb{R}$. Let $Q(z)dz^2$ be a quadratic
differential such that $\Gamma$ is a trajectory arc or an orthogonal
trajectory arc of $Q(z)dz^2$. Then let $D^{-}=\{\overline{z}:z\in
D\}$ be the reflection of $D$ along $\Gamma$. Define
\[Q^{-}(z)=\overline{Q(\overline{z})}\text{ for } z\in D^{-}.\]
Then since $\Gamma$ is a trajectory, we have
\[Q(z)=Q^{-}(z)\in\mathbb{R} \text{ for } z\in \Gamma.\]
Thus by defining
\begin{equation}Q^{*}(z)=\left\{\begin{array}{ll} Q(z) &\text{ for } z\in D, \\Q^{-}(z) &\text{ for } z\in D^{-},\\Q(z)=Q^{-}(z) &\text{ for } z\in \Gamma.\end{array}\right.\label{judah}\end{equation}
it is easy to see that $Q^{*}$ is meromorphic in $D\cup D^{-}$ and hence $Q^{*}(z)dz^2$ is a quadratic differential on $D\cup D^{-}$. Thus by the transformation law (and using Schwarz reflection), this shows that we can extend quadratic differentials across trajectory arcs or orthogonal trajectory arcs.

\ \\
We will use reflection to prove the following lemma:
\begin{lemma}\label{degbound}
Suppose $Q(z)dz^2$ is a K\"{u}hnau quadratic differential on $\mathbb{H}$. Then for any $z\in\partial\mathbb{H}$,
$\deg_{\mathbb{H},Q}(z)\in\mathbb{Z}$ implies that $Q(z)dz^2$ extends to a quadratic differential on a neighbourhood of $z$ and hence
\[\deg_{\mathbb{H},Q}(z)=\deg_{Q}(z).\]
\end{lemma}
\begin{proof}Firstly, if $z\in \Gamma_{j}$ for some $j=1,\ldots,n$. Then by definition $\deg_{\mathbb{H},Q}(z)=0$ and $Q(z)dz^2$ can be extended to a neighbourhood of $z$ by reflection. By definition, every $z\in\Gamma_{j}$ is an ordinary point of $Q(z)dz^2$ and hence  $\deg_{Q}(z)=\deg_{\mathbb{H},Q}(z)=0$.

Otherwise we write $z=z_{k}$ and suppose that a $\theta_{k-1}$-trajectory arc, $\Gamma_{k-1}$, and a $\theta_{k}$-trajectory arc, $\Gamma_{k}$, end at $z_{k}$.
Then, by definition, $\deg_{\mathbb{H},Q}(z_{k})\in\mathbb{Z}$ implies that $\theta_{k}-\theta_{k-1}$ is a multiple of $\pi/2$. Thus $\Gamma_{k-1}$ and $\Gamma_{k}$ are trajectory arcs or orthogonal trajectory arcs of $e^{-2i\theta_{k-1}}Q(z)dz^2$. Thus by reflection, $e^{-2i\theta_{k-1}}Q(z)dz^2$ extends to a neighbourhood of $z_{k}$. Hence, $Q(z)dz^2$ also extends to a neighbourhood of $z_{k}$.
\end{proof}

We can now prove Theorem \ref{QDSC}; but first, we explain briefly why we can view Theorem \ref{QDSC} as a generalized form of Schwarz-Christoffel mapping: Schwarz-Christoffel mapping is a method of computing the conformal map between the upper half-plane and a domain bounded by a polygon. See \cite{Neh82} for more details. If we have a
conformal map $f$ from $\mathbb{H}$ to some domain $D$  such that
 the sides of $D$ consist of $\theta$-trajectory arcs of the
quadratic differential $Q(w)dw^2$. Then $Q(w)dw^2$ is a K\"{u}hnau quadratic differential on $D$ and hence by Lemma \ref{gendeg}, $Q(f(z))f'(z)dz^2$ is a K\"{u}hnau quadratic differential on $\mathbb{H}$. Theorem \ref{QDSC} then implies that
\[Q(f(z))f'(z)^2=R\left(\prod_{j=1}^{n}(z-\zeta_{j})^{\lambda_{j}}\right).\]
This is precisely the Schwarz-Christoffel formula when $Q(z)\equiv 1$.

Also, we comment that the case when $Q(w)dw^2$ is either negative or positive on
$\mathbb{R}$ (i.e. the boundary of $\mathbb{R}$ consists only of
trajectory arcs and orthogonal trajectory arcs) is easy to prove: we
can use reflection to extend $Q(z)dz^2$ to a quadratic
differential on the Riemann sphere $\widehat{\mathbb{C}}$. Hence
$Q(z)$ must be rational since property (ii) in the definition of K\"{u}hnau quadratic differentials guarantees that $Q(z)$ does not have any essential singularities and so $Q(z)$ is rational (since the only meromorphic functions
on $\widehat{\mathbb{C}}$ are rational). This proves Theorem \ref{QDSC} for this case.
\begin{proof}[Proof of Theorem \ref{QDSC}]
Since $Q(z)dz^2$ is a K\"{u}hnau quadratic differential, we can find
\[z_{1}<\ldots<z_{m},\] and
\[\Gamma_{k}=\left\{\begin{array}{l} (z_{k-1},z_{k}) \text{ for }
 k=1,\ldots,m,\\ (z_{m},\infty) \text{ for } k=m+1, \\
 (-\infty,z_{0}) \text{ for } k=0, \end{array}\right.\] such that each
$\Gamma_{k}$ is a $\theta_{k}$-trajectory arc of $Q(z)dz^2$ for some $\theta_{k}\in[0,\pi)$. Let
\[\mathcal{T}=\left\{\Gamma_{1},\ldots,\Gamma_{m+1}\right\}.\]
Then take any $\Gamma\in \mathcal{T}$. Since $\Gamma$ is a
$\theta$-trajectory for some $\theta$, $\Gamma$ is a trajectory arc of $e^{-2i\theta}Q(z)dz^2$; hence by reflection, we can reflect the quadratic
differential $e^{-2i\theta}Q(z)dz^2$ across $\Gamma$ to get a
quadratic differential on $\mathbb{H}^{-}=\{\mathrm{Im}(z)<0\}$
which we call $\widetilde{Q}(z)dz^2$. Similarly, by rotating
$\widetilde{Q}(z)dz^2$, we can reflect it across another
$\Upsilon\in\mathcal{T}$ to get another quadratic differential
$Q^{*}(z)dz^2$ on $\mathbb{H}$. Since $Q^{*}$ is obtained from $Q$ by rotating twice, we have
\[Q^{*}(z)=e^{i\sigma}Q(z)\] for some $\sigma\in [0,2\pi)$. This
shows that
\[ \Psi(z)=\frac{Q'(z)}{Q(z)}=\frac{(Q^{*})'(z)}{Q^{*}(z)}\]
can be extended to a meromorphic function in $\mathbb{C}\setminus
\{z_{1},\ldots,z_{m}\}$. Then part (ii) of the definition of K\"{u}hnau quadratic differentials implies that all the finite singularities of
$\Psi(z)$ are simple poles otherwise $Q(z)$ would have an essential singularity which, by the great Picard theorem, contradicts part (ii) of the definition of K\"{u}hnau quadratic differentials. Thus we can write:
\[\Psi(z)=h(z) - \left(\sum_{j=1}^{n}
\frac{\lambda_{j}}{z-\zeta_{j}}\right),\]
where $\zeta_{j}\in \mathbb{C}$, and
$\lambda_{j}\in \mathbb{R}$, and $h(z)$ is
an entire function in $\mathbb{C}$ that does not does not vanish in $\mathbb{C}$.
This implies that
\[Q(z)=\exp\left[\int^{z}h(\zeta)d\zeta\right] \left(\prod_{j=1}^{n}(z-\zeta_{j})^{\lambda_{j}}\right)\left(\prod_{k=1}^{m}(z-z_{k})^{\nu_{k}}\right).\]
Moreover, the singularity at $\infty$ of
\[\exp\left[\int^{z}h(\zeta)d\zeta\right]\]
cannot be essential by part (ii) of the definition of K\"{u}hnau quadratic differentials (otherwise we would get a contradiction with the great Picard theorem as above). This implies that
\[\exp\left[\int^{z}h(\zeta)d\zeta\right]\]
is constant (since it has no zeroes or poles). Hence
\[Q(z)=R\left(\prod_{j=1}^{n}(z-\zeta_{j})^{\lambda_{j}}\right).\]
\end{proof}
If $\zeta_{j}\in\mathbb{C}\setminus\mathbb{R}$, then by definition, we must have
\[\lambda_{j}=\deg_{Q}(\zeta_{j}).\]
Moreover, if $\zeta_{j}\in\mathbb{R}$ and $\deg_{\mathbb{H},Q}(\zeta_{j})<\infty$ we also have
\[\nu_{k}=\deg_{\mathbb{H},Q}(z_{k}).\]
We will not prove this fact here but in the following corollary we will consider a special case. The general proof follows readily from it. We will prove the following corollary which is simply an application of Theorem \ref{QDSC} to domains slit by $\phi$-trajectory arcs:
\begin{corollary}\label{QDSCslit}
Suppose that $Q(z)dz^2$ is a K\"{u}hnau quadratic differential on
$\mathbb{H}$ such that there is a point $\xi_{0}\in\mathbb{R}$ with
$\deg_{\mathbb{H},Q}(\xi_{0})=N\in \{0,1,\ldots\}$; then we
can write
\begin{equation}Q(w)=R(w-\xi_{0})^N\left(\prod_{j=1}^{n}(w-a_{j})^{\alpha_{j}}
\right),\label{mundine}\end{equation}
 where $a_{j}\in \mathbb{C}$, $\alpha_{j}\in \mathbb{R}$, and also $R$ is some non-zero constant.
Let $\gamma:[0,T]\mapsto \overline{\mathbb{H}}$ be a simple curve such that $\gamma(0)=\xi_{0}$ and $\gamma(0,T)$ is a $\phi$-trajectory arc of
$Q(w)dw^{2}$ in $\mathbb{H}$ ($\phi\in[0,\pi)$) and $\zeta=\gamma(T)\in\mathbb{H}$ is an ordinary point of $Q(w)dw^2$ (i.e. $\deg_{Q}(\zeta)=0$). Suppose that $f$ maps $\mathbb{H}$ conformally onto $\mathbb{H}\setminus\gamma(0,T]$.
Then $f$ satisfies
\begin{equation}Q(f(z))f'(z)^2= R'(z-\xi)^2{(z-c^{-})^{\mu^{-}}(z-c^{+})^{\mu^{+}}}\prod_{j=1}^{n} (z-A_{j})^{\alpha_{j}}, \label{witter}\end{equation}
where $R'$ is some constant; $c^{-},c^{+}$ are the two preimages of
$\xi_{0}$ under $f$ satisfying $c^{-}<c^{+}$;
$A_{j}$ is the preimage of $a_{j}$ under $f$; and $\xi$ is the
preimage of $\zeta$; and
\[\mu^{\pm}=\deg_{\mathbb{H}\setminus\gamma(0,s],Q}(f(c^{\pm})).\]
\end{corollary}
\begin{proof}Theorem \ref{QDSC} and Lemma \ref{degbound} imply that $Q(w)$ can be written as (\ref{mundine}).
Then Lemma \ref{gendeg} implies that $\widehat{Q}(z)=Q(f(z))f'(z)^2$ is a K\"{u}hnau quadratic
differential. So by Theorem \ref{QDSC}, we only need to look at the singularities of $\widehat{Q}(z)dz^2$.

Now, by Schwarz reflection, $f$ extends to a conformal map on $\mathbb{C}\setminus\{\xi,c^{-},c^{+}\}$. Thus, by the conformal invariance of trajectories, this implies that
\[Q_{f}(z)= R'(z-\xi)^M{(z-c^{-})^{\mu^{-}}(z-c^{+})^{\mu^{+}}}\prod_{j=1}^{n} (z-A_{j})^{\alpha_{j}}.\]
Then by Lemma \ref{loc}, there are exactly two $\phi$-trajectory arcs of $Q(z)dz^2$ ending at $\zeta=\gamma(T)$ of which $\gamma(0,T)$ is one of them. So by the conformal invariance of trajectories, there is one $\phi$-trajectory arcs of $\widehat{Q}(z)dz^2$ ending at $\xi$ that is contained in $\mathbb{H}$. Hence, by definition, $\deg_{\widehat{Q},\mathbb{H}}(\xi)=2$. Using Lemma \ref{degbound}, this implies that $\deg_{Q}(\xi)=2$ i.e. $M=2$. Thus we only need to determine $\mu^{-}$ and $\mu^{+}$.

Note that since $\xi_{0}$ has degree $N$ with respect to $Q(z)dz^2$, we can determine, using Lemma \ref{loc}, that the angle between $\gamma(0,T)$ and $f((\xi,c^{-}) )$ at $\xi_{0}$ is
\[\pi\psi^{-}=\pi \frac{\deg_{H,Q}(f(c^{-}))+2}{N+2},\] and similarly, the angle between
$\gamma(0,T)$ and $f( (\xi,c^{-}) )$ at $\xi_{0}$ is
\[\pi\psi^{+}= \pi\frac{\deg_{H,Q}(f(c^{+}))+2}{N+2}.\]
Hence, by Schwarz reflection, the function
\[F(z)=(f(z)-\xi_{0})^{1/\psi^{-}}\]
extends to a conformal mapping on a neighbourhood of $c^{-}$. Thus in a neighbourhood of $z=c^{-}$, we can write
\begin{equation}f(z)=\xi_{0}+(z-c^{-})^{\psi^{-}}h(z)^{\psi^{-}},\label{chagaev}\end{equation} where $h$ is
analytic in a neighbourhood of $z=c^{-}$ with $h(c^{-})\neq0$. Now
\[\frac{Q_{f}'(z)}{Q_{f}(z)}=\frac{Q'(f(z))f'(z)^2}{Q(f(z))}+2\frac{f''(z)}{f'(z)}.\]
The residue at $z=c^{-}$ of the left-hand side of the equation is
$\mu^{-}$, and we can use (\ref{chagaev}) to determine the residue
at $z=c^{-}$ of the right-hand side. Thus we get
\[\mu^{-}=\deg_{H,Q}(f(c^{-})).\] We apply the same method to $c^{+}$ to
get $\mu^{+}$.
\end{proof}

\section{Domains slit by $\theta$-trajectory arcs}
Let $Q(w)dw^2$ be a K\"{u}hnau quadratic differential on $\mathbb{H}$ with $\deg_{\mathbb{H},Q}(\xi_{0})=N\in \{0,1,\ldots\}$ for some $\xi_{0}\in\mathbb{R}$. Then by Theorem \ref{QDSC} and Lemma \ref{degbound},
\[Q(w)=(w-\xi_{0})^N\left(\prod_{j=1}^{n}(w-a_{j})^{\alpha_{j}}\right),\]
where $\alpha_{j}\in\mathbb{R}$ and $a_{j}\in\mathbb{C}$. Now
suppose that $\gamma:[0,T)\mapsto \overline{\mathbb{H}}$ is a simple curve such that $\gamma(0)=\xi_{0}$ and $\gamma(0,T)$ is a  $\phi$-trajectory of $Q(z)dz^2$ in $\mathbb{H}$ ($\phi=[0,\pi)$) that is parameterized by
half-plane capacity. As mentioned
in the introduction, there exists conformal maps $f_{t}:\mathbb{H}\mapsto
H_{t}=\mathbb{H}\setminus\gamma(0,t]$ satisfying the hydrodynamic
normalization. Then by restricting $Q(w)dw^2$ to a quadratic
differential on $H_{t}$ we can induce via $f_{t}$ and (\ref{khan}),
a quadratic differential on $\mathbb{H}$:
\begin{equation}Q_{t}(z)dz^2=Q(f_{t}(z))f'_{t}(z)^2dz^2.\label{calzaghe}\end{equation}
We now use Corollary \ref{QDSCslit} and (\ref{calzaghe}) to prove Theorem \ref{QDdriving}.
\begin{proof}[Proof of Theorem \ref{QDdriving}]
Note that by Schwarz reflection, each $f_{t}$ can be extended to a conformal map on $\widehat{\mathbb{C}}\setminus\{C^{-}(t),C^{+}(t),\xi(t)\}$.
Then since $f_{t}(z)$ satisfies the hydrodynamic normalization, this implies that
\[f'_{t}(z)^{2}=1+O\left(\frac{1}{z^{2}}\right) \text{ as } z\rightarrow\infty.\]
So by (\ref{calzaghe}),
\begin{equation}\frac{Q_{t}(z)}{Q(f_{t}(z))}=1+O\left(\frac{1}{z^{2}}\right) \text{ as } z\rightarrow\infty.\label{tarver}\end{equation}
If we let $\zeta=1/z$, then we get
\begin{equation}\frac{Q_{t}(1/\zeta)}{Q(f_{t}(1/\zeta))}=1+O\left(\zeta^2\right) \text{ as } \zeta\rightarrow0.\label{kessler}\end{equation}
Since $f_{t}$ is analytic in a neighbourhood of infinity, (\ref{kessler}) is a Taylor series expansion and hence we can look at the Taylor series coefficients, in particular:
\[\int_{C(0,\epsilon)} \frac{f'(1/\zeta)^2}{\zeta^2} d\zeta =\int_{C(0,\epsilon)} \frac{Q_{t}(1/\zeta)}{\zeta^{2}Q(f_{t}(1/\zeta))} d\zeta =0\]
for small enough $\epsilon >0$ where $C(0,\epsilon)$ is the
anticlockwise contour about the circle with centre at zero and
radius $\epsilon>0$.
Then by Theorem \ref{QDSC} and Corollary \ref{QDSCslit}, we can write
\[Q(w)=R(w-\xi_{0})^N\left(\prod_{j=1}^{n}(w-a_{j})^{\alpha_{j}}
\right),\]
and
\[Q_{t}(z)=R'(z-\xi)^2{(z-C^{-}(t))^{\mu^{-}}(z-C^{+}(t))^{\mu^{+}}}\prod_{j=1}^{n} (z-A_{j}(t))^{\alpha_{j}}.\]
Hence by the residue theorem (since
$f_{t}(1/\zeta)=1/\zeta+\cdots$ as $\zeta\rightarrow 0$), this
implies that
\small\[\left[2\xi(t) + \mu^{-}C^{-}(t) +\mu^{+}C^{+}(t) + \left(\sum_{k=1}^{n} \alpha_{k}A_{k}(t)\right)\right] \\- \left[N\xi_{0}+  \left(\sum_{k=1}^{n} \alpha_{k}a_{k}\right) - \left(\sum_{l=1}^{m} \beta_{l}b_{l}\right)  \right]=0.\]\normalsize
This implies (\ref{toney}). To get (\ref{wright}), note that $f_{t}$ satisfies
the chordal Loewner differential equation (\ref{hatton}) and hence if we let
 $G_{t}=f_{t}^{-1}\circ f_{s}$ for some $s\in (0,T)$ fixed and $t>s$, then the chain rule implies that $G_{t}$ satisfies the differential equation
\[\frac{\partial G_{t}}{\partial t}(z)= \frac{2}{G_{t}(z)-\xi(t)}.\]
Then for some $s$ sufficiently close to $t$, we can write each $A_{j}(t)=G_{t}(w_{j})$ $w_{j}\in \mathbb{C}$ for all $j=1,\ldots,n$. Thus
\[\dot{A}_{j}(t)=\frac{2}{A_{j}(t)-\xi(t)}.\]
Similarly, we get
\[ \dot{C}^{\pm}(t)=\frac{2}{C^{\pm}(t)-\xi(t)}.\]
Hence we get (\ref{wright}) from differentiating (\ref{toney}).
\end{proof}
\subsection*{An extension:}\ \\
We can extend Theorem \ref{QDdriving} to the case when $\gamma$ is made up of
different $\theta_{k}$-trajectory arcs of some quadratic
differential $Q(z)dz^2$: Let $\gamma:(0,T]\mapsto\mathbb{H}$ be a
curve with $\gamma(0)\in \mathbb{R}$ such that there is a partition
\[\{0=t_{0}<t_{1}<\cdots<t_{r}=T\}\] such that
$\gamma(t_{k-1},t_{k})$ is a $\theta_{k}$-trajectory arc of
$Q(z)dz^2$ and $\gamma(t_{k})$ is an ordinary point of $Q(z)dz^2$ for $k=1,\ldots, r$. Then we can find the driving function $\xi(t)$ of
$\gamma$ by applying Theorem \ref{QDdriving} to the $\theta_{1}$-trajectory arc
$\gamma(0,t_{1})$ to get a driving function $\xi_{1}(t)$, and
applying Theorem \ref{QDdriving} inductively to each
$f_{t_{k}}^{-1}(\gamma(t_{k},t_{k+1}))$ (which is a
$\theta_{k+1}$-trajectory arc of the quadratic differential
$Q_{t_{k}}(z)dz^2=Q(f_{t_{k}}(z))f'_{t_{k}}(z)^2dz^2$) to get
$\xi_{k}(t)$. Then
\[\xi(t)= \xi_{k}(t) \text{ for } t\in [t_{k-1},t_{k}).\]
We also have the following corollary:
\begin{corollary}\label{smoothxi}
Suppose that $Q(w)dw^2$ and $\gamma$ are
as defined in Theorem \ref{QDdriving}. Then the driving function $\xi$ and
$A_{j},C^{-}, C^{+}$ as defined in Theorem \ref{QDdriving} are in
$C^{\infty}(0,T)$. Moreover, we can write any derivative  of $\xi,C^{-},
C^{+},A_{j}$ explicitly in terms of $\xi,C^{-}, C^{+},A_{j}$ and the
exponents $\mu^{-},\mu^{+},\alpha_{j}$.
\end{corollary}
\begin{proof} Recall that, in the proof of Theorem \ref{QDdriving}, we had the formulae
\[\dot{A}_{j}(t)=\frac{2}{A_{j}(t)-\xi(t)}, \dot{C}^{\pm}(t)=\frac{2}{C^{\pm}(t)-\xi(t)}.\]
This implies that each term in (\ref{wright}) is differentiable so
we can write the second derivative of $\xi$ in terms of
$\xi(t),A_{j}(t),C^{\pm}(t)$. This in turn implies that we can write
the third derivative of $\xi$ in terms of
$\xi(t),A_{j}(t),C^{\pm}(t)$ and the exponents. Continuing
inductively, we have showed that every derivative of $\xi$ exists
and can be expressed in terms of $\xi(t),A_{j}(t),C^{\pm}(t)$ and
the exponents. Note that each derivative of $\xi$ is finite for
$t\in(0,T)$ since
\[ |A_{j}(t)-\xi(t)|,|C^{\pm}(t)-\xi(t)|>0.\]
Then $\xi$ is smooth implies that $A_{j}(t),C^{\pm}(t)$ are also smooth.
\end{proof}
Theorem \ref{gammadot} then follows from Corollary \ref{QDSCslit} and Theorem \ref{QDdriving}:
\begin{proof}[Proof of Theorem \ref{gammadot}]\ \\
First note that, by the definition of $\theta$-trajectory arcs, $\dot{\gamma}$ always exists and is never 0. Also by Corollary \ref{QDSCslit}, $\Phi_{t}(\xi(t))\neq 0,\infty$; thus the right hand side of (\ref{erdei}) always exists since, by definition, $\gamma$ avoids poles and zeroes of  $Q(w)dw^2$.

Recall that $f_{t}(\xi(t))=\gamma(t)$, this implies that
\[\dot{\gamma}(t) = \dot{f}_{t}(\xi(t)) + f_{t}'(\xi(t))\dot{\xi}(t).\]
Then combining the Loewner differential equation (\ref{hatton}) with (\ref{calzaghe}) we have
\[\dot{f}_{t}(z)=-\frac{2}{z-\xi(t)}\sqrt{\frac{Q_{t}(z)}{Q(f_{t}(z))}}=-2\sqrt{\frac{\Phi_{t}(z)}{Q(f_{t}(z))}}\]
\[\Rightarrow \dot{f}_{t}(\xi(t))=-2\sqrt{\frac{\Phi_{t}(\xi(t))}{Q(\gamma(t))}}.\]
Note that $\Phi_{t}(\xi(t))\neq 0,\infty$ since, by Corollary \ref{QDSCslit}, $Q_{t}(z)$ has a double zero at $\xi(t)$.

Thus we have
\[\dot{\gamma}(t) = -2\sqrt{\frac{\Phi_{t}(\xi(t))}{Q(\gamma(t))}} + \sqrt{\frac{Q_{t}(\xi(t))}{Q(\gamma(t))}}\dot{\xi}(t)=-2\sqrt{\frac{\Phi_{t}(\xi(t))}{Q(\gamma(t))}},\] since Theorem \ref{QDdriving} implies that $\dot{\xi}$ is finite for all $t\in(0,T)$ and Corollary \ref{QDSCslit} implies that $Q_{t}(\xi(t))=0$.
\end{proof}
\section{Applying Theorem \ref{QDdriving}}
In practice, understanding $\xi(t)$ via (\ref{toney}) is not
possible: it is difficult to calculate the positions of the zeroes
and poles of $Q_{t}$ because the information we have on them is all
relative to $\xi(t)$ (which we are trying to find). On the other
hand, (\ref{wright}) is more useful in applications. In this
section, we will demonstrate how we can use (\ref{wright}) to
calculate numerically the driving function of a given slit that consists of
$\theta_{k}$-trajectory arcs of a given quadratic differential. The method is basically a modified version of Euler's method.

Firstly, for any smooth function $h$ on (0,T), Taylor's theorem implies that for all $M=1,2,\ldots$,
\begin{equation}
\left|h\left(t+\frac{1}{K}\right)-\left( h(t) + \sum_{m=1}^{M-1}
\frac{1}{m!K^m}\frac{d^{m}h}{dt^{m}}(t) \right)\right| \leq
\frac{1}{M!K^{M}}\sup_{s\in\left(t,t+\frac{1}{K}\right)}
\left|\frac{d^{M}h}{dt^{M}}(s)\right| \label{castillo}\end{equation}
for $t,t+1/K\in(0,T)$. We will apply (\ref{castillo}) to the
functions $\xi,A_{k}$ and $C^{\pm}$ (as defined in Theorem \ref{QDdriving})
noting that, by Corollary \ref{smoothxi}, they are smooth and all of their derivatives can be expressed in terms of $\xi(t), A_{k}(t)$ and $C^{\pm}(t)$. Thus
if we know $\xi(s), A_{k}(s)$ and $C^{\pm}(s)$ we can
use (\ref{castillo}) to obtain an approximate formula for
$\xi(s+K^{-1}), A_{k}(s+K^{-1})$ and $C^{\pm}(s+K^{-1})$ (choosing $K$
to be small and/or $M$ to be large so that the right-hand-side of (\ref{castillo}) is small); then we can apply
(\ref{castillo}) to $\xi(s+K^{-1}), A_{k}(s+K^{-1})$ and
$C^{\pm}(s+K^{-1})$ to find $\xi(s+2K^{-1}), A_{k}(s+2K^{-1})$ and
$C^{\pm}(s+2K^{-1})$. Continuing like this, we obtain an approximation
of $\xi$ at the points $\{s+nK^{-1}\}$.

So clearly what we need to do now is find the starting values
$\xi(s), A_{k}(s)$ and $C^{\pm}(s)$ so we can apply the above
method. But because $\xi$ is not differentiable at $0$, we cannot
use the formula (\ref{castillo}) with $t=0$. The way around this is
to note that if $\deg_{\mathbb{H},Q}(\xi_{0})=N\in \{0,1,2,\ldots\}$
then since we know $\deg_{\mathbb{H},Q_{t}}(C^{+}(t))$, we
can calculate the angle that the trajectory makes with the line
$[\xi_{0},\infty)$ (as in the proof of Corollary \ref{QDSCslit}). Then we find
that the angle is $\pi\psi$ where:
\[ \pi\psi=\pi\left(\frac{2\deg_{\mathbb{H},Q_{t}}(C^{+}(t))+2}{N+2}\right).\]
So if we choose $s$ small enough, we have
\[f_{s}\approx F_{s}^{\psi,\xi_{0}},\]
where $F_{s}^{\psi,\xi_{0}}$ is the conformal map that maps
$\mathbb{H}$ conformally onto $H^{\psi}_{s}$ that is
hydrodynamically normalized where $H^{\psi}_{s}$ is the upper
half-plane slit by the straight line starting at $\xi_{0}$ making an
angle $\pi\psi$ with $[\xi_{0},\infty)$, with half-plane
capacity $2s$. Then we also have
\[A_{k}(s)\approx(F_{s}^{\psi^{+},\xi_{0}})^{-1}(a_{k}),\]
and also, $C^{-}(s),C^{+}(s)$ are approximately the two preimages of
$\xi_{0}$ under $F_{s}^{\psi^{+},\xi_{0}}$. Then we can use
(\ref{toney}) to calculate $\xi(s)$ approximately. We can then
plug this information into (\ref{castillo}) as described above.

Note that $F_{t}^{\psi,x}$ can be found using the fact that
\begin{equation} F_{\lambda t}^{p,0}(z)=\left(z-(1-2p)\sqrt{t}-p\sqrt{t}\right)^p\left(z-(1-2p)\sqrt{t}+(1-p)\sqrt{t}\right)^{1-p}\label{harrison}\end{equation}
for some $\lambda$. Then we reparameterize this formula to remove the $\lambda$
and translate the point 0 to $\xi_{0}$. Unfortunately, inverting this
function cannot be done explicitly but it can be done numerically
very efficiently using Newton's method. Alternatively, by selecting a small $s$, we can assume that
\[A_{k}(s)\approx a_{k}\]
for all $k$. Then we note that the 2 preimages of $\xi_{0}$ under $F_{s}^{\psi,x}$ can be determined explicitly (see \cite{MR06}). This obviates the need to numerically invert $F_{t}^{\psi,x}$.
\ \\

Another difficulty is that, in general, given a slit, we cannot
parameterize it by half-plane capacity so it would be difficult, for
example, to know at which $t$ one should stop. Most formulae for
calculating half-plane capacity of some compact set $K$ rely on
knowing the conformal map $f_{K}$ of $\mathbb{H}$ onto
$\mathbb{H}\setminus K$ (normalized hydrodynamically). One possibility would be to use the probabilistic definitions of half-plane
capacity given in \cite{Law05}. We will use the fact that Theorem \ref{gammadot} and Corollary \ref{smoothxi} imply that we can give all derivatives of
$\gamma(t)$ in terms of $\xi(t),A_{k}(t),C^{-}(t),C^{+}(t)$ and the
exponents $\mu^{-},\mu^{+}, \alpha_{k}$ so if we know these, we can
also use (\ref{castillo}) to approximate $\gamma$. This in turn
allows us to calculate the length of the slit $\gamma$. Thus if we know
beforehand length of our slit, we can calculate at what value of $t$
we stop.

We now have everything we need in order to use (\ref{castillo})
to calculate the driving function numerically of any slit that is
made up of $\theta_{k}$ trajectory arcs of a quadratic differential
$Q(w)dw^2$. We will demonstrate how this is done in the following
example:

\subsection*{An example.}
Suppose that $\gamma:(0,T)\rightarrow \mathbb{H}$ is a piecewise
linear arc parameterized by half-plane capacity that satisfies:
\begin{itemize}
 \item $\gamma(0)=0$.
\item From $t=0$ to $t=t_{1}$, $\gamma$ is the straight line arc from $0$ to $i$; call this $\Gamma_{1}$.
\item From $t=t_{1}$ to $t=t_{2}$, $\gamma$ is the straight line arc from $i$ to $2+i$; call this $\Gamma_{2}$.
\item From $t=t_{2}$ to $t=t_{3}=T$, $\gamma$ is the straight line arc from $2+i$ to $2+2i$; call this $\Gamma_{3}$.
\end{itemize}
First note that $\gamma$ is made up of alternating $(\pi/2)$- and $0$-trajectory arcs of the quadratic differential $1dw^2$ in $\mathbb{H}$ and hence we can use Theorem \ref{QDdriving} (or more specifically the extension of Theorem \ref{QDdriving} detailed in Section 3) to calculate $\dot{\xi}$. As mentioned previously, there is no easy way to know beforehand what $t_{1},\ldots,t_{3}$ are.
For simplicity, we will only use $M=1$ in (\ref{castillo}) i.e.
\[f\left(t+\frac{1}{K}\right)\approx f(t)+ \frac{\dot{f}(t)}{K},\]
and fix a large $K$.
Obviously $\Gamma_{1}$ forms a right angle with real line; so we can use (\ref{harrison}) to determine the function \[f_{t_{1}}=F_{t_{1}}^{1/2,0}(z)=\sqrt{z^2-4t_{1}}.\]

\pagebreak
\begin{figure}[hp]
 \begin{center}
\scalebox{0.5}{\includegraphics{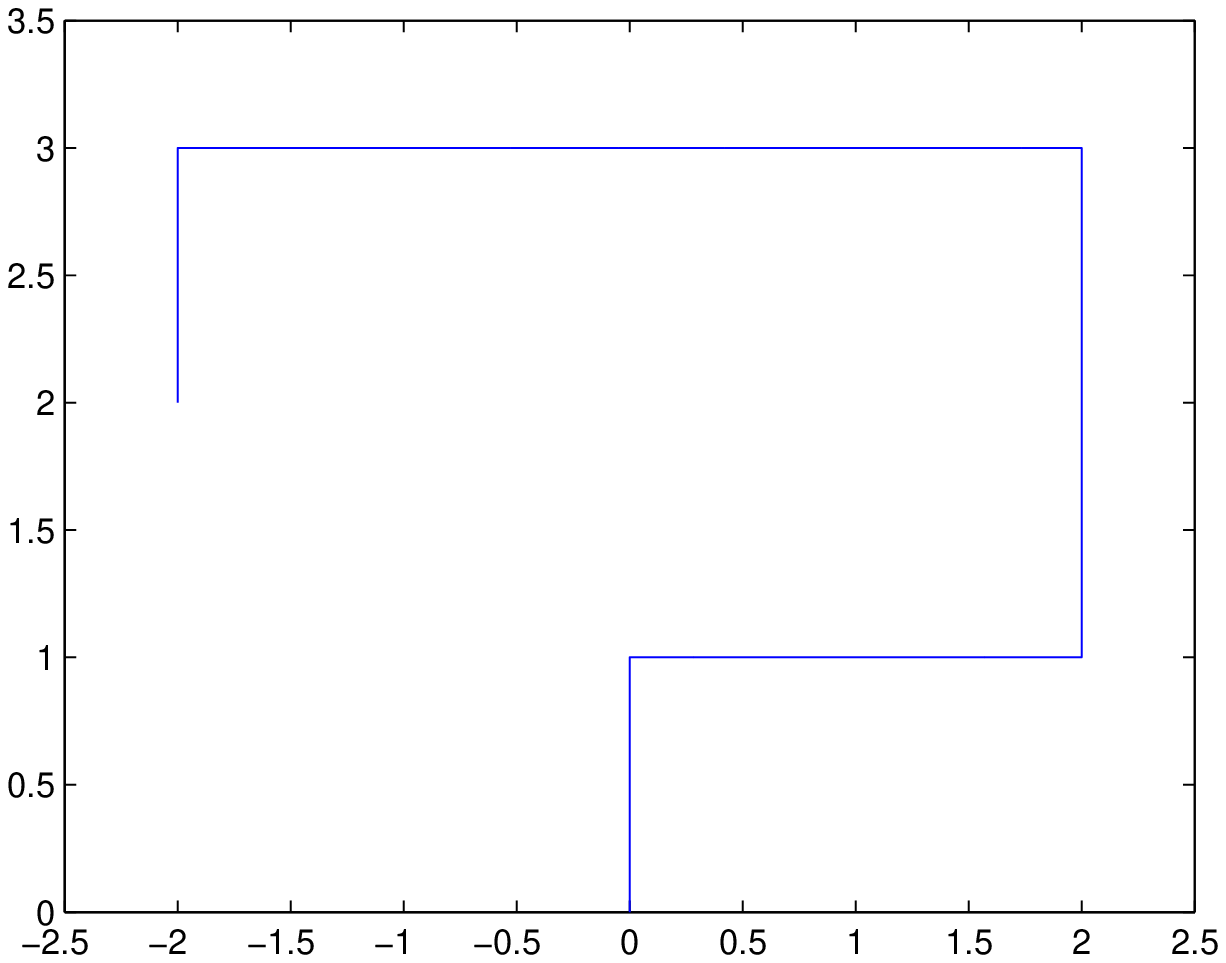}\includegraphics{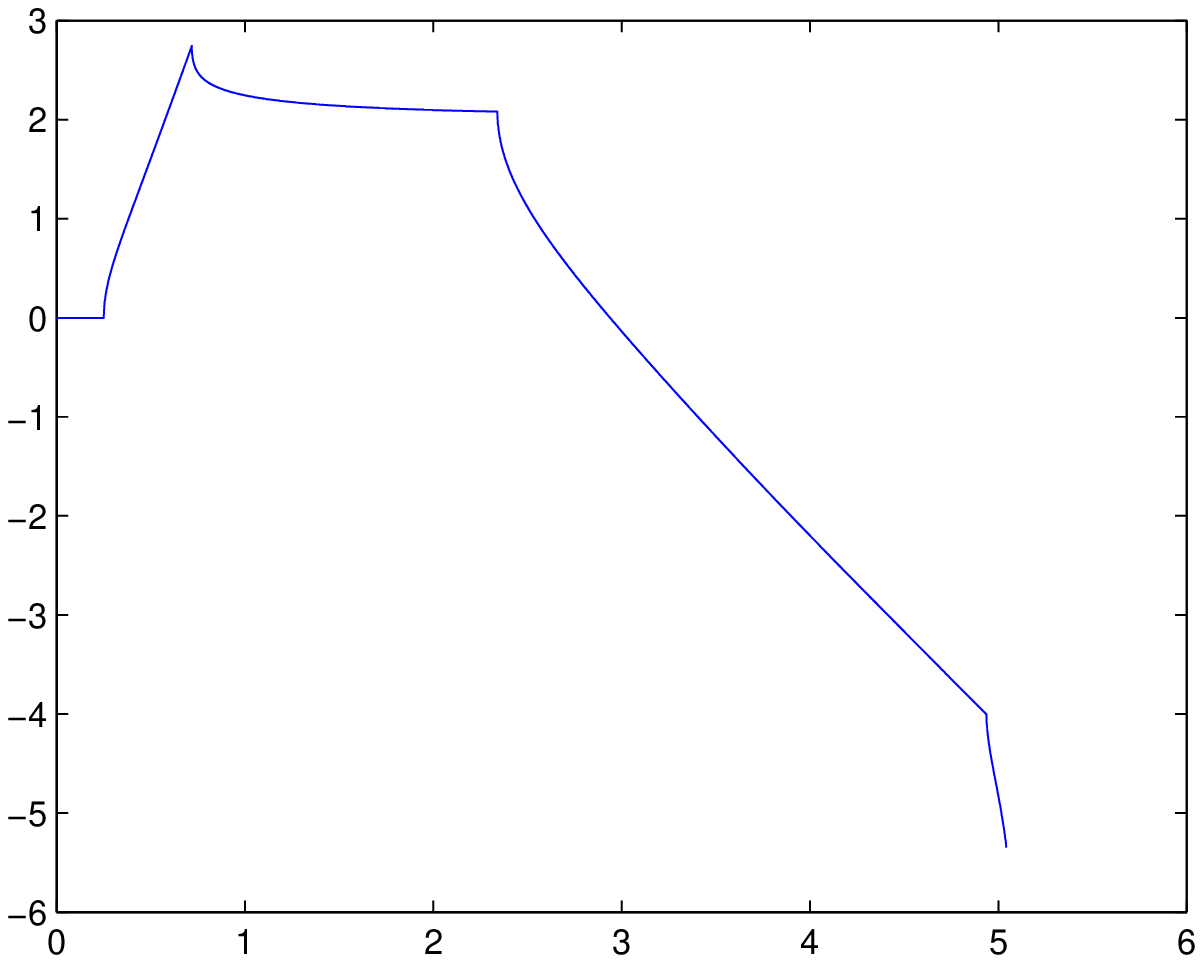}}
 \end{center}
\caption{The example path in the upper half-plane (left) and a plot
of its driving function on the $y$-axis against time on the $x$-axis
(right).} \label{fig2}
\end{figure}
It is easy to see that in this case, $t_{1}=1/4$ and $\xi$ is constantly 0 for $t\in(0,t_{1}]$. This induces the quadratic differential using (\ref{calzaghe}):
\[Q_{t_{1}}(z)dz^2=\frac{z^2dz^2}{(z+1)(z-1)}.\]
Hence, we let $A_{1}(t_{1})=-1$, $A_{2}(t_{1})=1$. Also $f_{t_{1}}^{-1}(\gamma_{2})$ is a $0$-trajectory arc of $Q_{1}(z)dz^2$ starting from $\xi(t_{1})=0$ on $\mathbb{R}$ (by the conformal invariance of trajectories). Now note that $f_{t_{1}}^{-1}(\gamma_{2})$ makes an angle of $\pi/4$ with the positive real axis. and hence
\[f_{t_{1}+K^{-1}}\approx F_{K^{-1}}^{1/4,0}(z)\]
since $K$ is large. We can then use Newton's method to find the preimages under the above approximation of $f_{t_{1}+K^{-1}}$
 of the points $A_{1}(t_{1}), A_{2}(t_{1})$ and the
2 preimages of zero to get the points $A_{1}(t_{1}+K^{-1}), A_{2}(t_{1}+K^{-1}),C^{-}(t_{1}+K^{-1}),C^{+}(t_{1}+K^{-1})$ and hence, using (\ref{toney}), we can find $\xi(t_{1}+K^{-1})$.
Then inserting this into (\ref{castillo}), as detailed above we can also find $\xi(t_{1}+nK^{-1})$ and $A_{1}(t_{1}+nK^{-1}), A_{2}(t_{1}+nK^{-1}),C^{-}(t_{1}+nK^{-1}),C^{+}(t_{1}+nK^{-1})$; also, by Theorem \ref{gammadot}, we can find $|\dot{\gamma}(t_{1}+nK^{-1})|$  if we let
\[t_{2}(K)= \inf\left\{n: \sum_{j=1}^{n} \frac{1}{K}|\dot{\gamma}(t_{1}+nK^{-1})|> (\text{length of } \Gamma_{2})=2\right\},\]
then $t_{2}(K)\approx t_{2}$ for $K$ large. So we just assume that
$t_{2}=t_{2}(K)$. Let $A_{3}(t_{2})=C^{-}(t_{2})$ and
$A_{4}(t_{2})=C^{+}(t_{2})$. Hence by (\ref{calzaghe}),
\[Q_{t_{2}}(z)dz^2= \frac{(z-\xi(t_{2}))^2(z-A_{3}(t_{2}))dz^2}{(z- A_{1}(t_{2}))(z- A_{2}(t_{2}))(z- A_{4}(t_{2}))}.\]
Then, by the conformal invariance of trajectories, $f_{t_{2}}^{-1}(\Gamma_{3})$ is a $\pi/2$-trajectory of
$Q_{t_{2}}(z)dz^2$ and also, $f_{t_{2}}^{-1}(\Gamma_{3})$, makes an
angle $3\pi/4$ with $(\xi(t_{2}),\infty)$ and so
\[f_{t_{2}+K^{-1}}\approx F_{K^{-1}}^{3/4,0}(z).\]
Then, as before, we can use Newton's method to find the preimages
under the above approximation of $f_{t_{2}+K^{-1}}$
 of the points $A_{1}(t_{2}),\ldots, A_{4}(t_{2})$ and the
2 preimages of $\xi(t_{2})$ to get the points $A_{1}(t_{2}+K^{-1}),\ldots,
A_{4}(t_{2}+K^{-1}),C^{-}(t_{2}+K^{-1}),C^{+}(t_{2}+K^{-1})$ and
hence use (\ref{toney}) to get $\xi(t_{2}+K^{-1})$. We insert these
into the formula iteratively to get $\xi(t_{2}+nK^{-1})$ and
$A_{1}(t_{2}+nK^{-1}),\ldots,
A_{4}(t_{2}+nK^{-1}),C^{-}(t_{2}+nK^{-1}),C^{+}(t_{2}+nK^{-1})$
until $t_{2}+nK^{-1}\approx T$. Thus the end result is that we found
the driving function of the first 3 steps of the slit given in
Figure \ref{fig2}. Of course, our calculation of $\xi$ will be more accurate by taking larger $K$.
\ \\

 For example, we
can use the above method to calculate the driving function of any
path on the square/triangle/hexagonal lattice on $\mathbb{H}$
starting from some point in $\mathbb{R}$. In fact we can calculate
the driving function of a path on the square/triangle/hexagonal
lattice in any polygon $D$ by mapping  the half-plane conformally
onto $D$ and pulling back the quadratic differential $1 dw^2$ on $D$
to $Q(z)dz^2$ on $\mathbb{H}$ using the transformation law. Also
note that, in general, any curve $\gamma$ can be approximated by a
curve $\gamma_{\delta}$ which lies on the square lattice
$\delta\mathbb{Z}^{2}$. Then it can be shown that
\[\xi_{\delta}\rightarrow \xi \text{ uniformly as } \delta\searrow 0,\]
where $\xi_{\delta}$ is the driving function of $\gamma_{\delta}$
and $\xi$ is the driving function of $\gamma$ hence, we can use the
above method to calculate $\xi_{\delta}$ then take the limit as
$\delta\searrow0$ to obtain $\xi$.

Another point to note is that using the above method, we do not need to know before hand what the trajectory arc of the given quadratic differential looks like; so for arbitrary K\"{u}hnau quadratic differentials, we can use this method to plot the trajectories starting at the boundary.

We end this section by looking at what happens when the slit approaches the boundary:

\begin{proposition}
Suppose that $\gamma:[0,T)\mapsto \overline{\mathbb{H}}$ is a simple curve such that $\gamma(0)\in\mathbb{R}$ and $\gamma(0,T)$ is a $\theta$-trajectory arc of some quadratic differential $Q(z)dz^2$. Then let $\xi$ be the driving function of $\gamma$. If
\[\lim_{t\uparrow T}\gamma(t)\in \mathbb{R}\cup \gamma(0,T),\]
i.e. $\gamma$ makes a loop at time $T$. Then
\[\left|\frac{d^{n}\xi}{d^{n}t}(t)\right|\rightarrow \infty \text{
as } t\nearrow T\] for all $n=0,1,\ldots$.
\end{proposition}
\begin{proof}\ \\
For $t\in(0,T)$, we define \[\Gamma(t)=\{\gamma(s):s\in(t,T)\}.\]
Then $\Gamma_{t}$ is a $\theta$-trajectory arc in
$H_{t}=\mathbb{H}\setminus\gamma(0,t]$ of $Q(w)dw^2$ and it is also
a crosscut in $H_{t}$ (see \cite{Pom92}). Then by the conformal
invariance of $\theta$-trajectories,
$f_{t}^{-1}(\Gamma_{t})\subset\mathbb{H}$ is a $\theta$-trajectory
arc of $Q_{t}(z)dz^2$. Moreover, $f_{t}^{-1}(\Gamma_{t})$ is a
crosscut of $\mathbb{H}$ with one end point at $\xi(t)$ and the other
end point in $\mathbb{R}$ such that either $C^{+}(t)$ or $C^{-}(t)$
is contained in the closure of the bounded component of
$\mathbb{H}\setminus f_{t}^{-1}(\Gamma_{t})$. Without loss of
generality, assume it is $C^{+}(t)$. Then since
$\mathrm{diam}(f_{t}^{-1}(\Gamma_{t}))\rightarrow 0$ as $t\nearrow
T$, we must have $\xi(t)=C^{+}(T)$ and hence by (\ref{wright}),
$\dot{\xi}(t)\rightarrow \infty$ as $t\nearrow T$. Similarly, we
differentiate (\ref{wright}) as mentioned in Corollary \ref{smoothxi} to obtain
the result for higher order derivatives.
\end{proof}
This means that as $\gamma$ gets closer and closer to making a loop,
the approximation by (\ref{castillo}) stops working no matter what
$M$ we choose. This phenomenon can be observed in Figure \ref{fig2},
as we turn the last corner in $\gamma$, we can see that $\xi$
decreases faster even though the slit is not yet that
close to the boundary.

\section{Generalizing Theorem \ref{QDdriving}}
\subsection{Multiple slits}
Suppose that $\gamma_{k}:[0,T)\rightarrow \overline{\mathbb{H}}$ for $k=1,\ldots,N$
are disjoint simple curves such that
$\gamma_{k}(0)\in\mathbb{R}$ and $\gamma_{k}(0,T)\subset\mathbb{H}$. By the Riemann mapping theorem, there exists
unique $f_{t}$ that map $\mathbb{H}$ conformally onto
$H_{t}=\mathbb{H}\setminus \bigcup_{k=1}^{N}\gamma_{k}(0,t]$ that
satisfies the hydrodynamic normalization. We can reparameterize such
that \[\bigcup_{k=1}^{N}\gamma_{k}(0,T]\] has half-plane capacity
$2t$. Then $f_{t}$ satisfies
\begin{equation}\dot{f}_{t}(z)=-2f_{t}'(z)\left(\sum_{k=1}^{N}\frac{b_{k}(t)}{z-\xi_{k}(t)}\right),\label{urango}\end{equation}
where
\[\sum_{k=1}^{N} b_{k}(t) = 1,\]
and $\xi_{k}(t)=f_{t}^{-1}(\gamma_{k}(t))$. See \cite{KNK04} for
more details.

\begin{theorem}
Suppose that $Q(w)dw^2$ is a K\"{u}hnau quadratic differential on
$\mathbb{H}$ such that the points $\xi_{k}(0)\in\mathbb{R}$ satisfy
\[\deg_{\mathbb{H},Q}(\xi_{k}(0))=\beta_{k}\in \{0,1,2,\ldots\}\]
for all $k$. Then we can write
\[Q(w)=\left(\prod_{k=1}^{N}(w-\xi_{k}(0))^{\beta_{k}}\right)\left(\prod_{j=1}^{n}(w-a_{j})^{\alpha_{j}}   \right)  \]
with $a_{j}\in\mathbb{C}$ and $\alpha_{j}\in\mathbb{R}$. Then suppose that $\gamma_{k}:[0,T)\rightarrow \overline{\mathbb{H}}$ for $k=1,\ldots,N$
are disjoint simple curves such that
$\gamma_{k}(0)\in\mathbb{R}$ and $\gamma_{k}(0,T)\subset\mathbb{H}$ and are parameterized as above. Then
\begin{equation}
2\sum_{k=1}^{N}\xi_{k}=-\left(\sum_{k=1}^{N}(\mu_{k}^{-}C_{k}^{-}(t)+\mu_{k}^{+}C_{k}^{+}(t))\right)-\left(\sum_{j=1}^{n}\alpha_{j}A_{j}(t)\right)+\Sigma_{0},
\label{rahman}\end{equation} and
\small\begin{equation}\dot{\xi}_{l}(t)= \left(\sum_{k=1,k\neq l}^{N}\frac{b_{k}(t)}{\xi_{l}(t)-\xi_{k}(t)}\right) -\frac{1}{2}\left(\sum_{k=1}^{N}\frac{\mu_{k}^{-}b_{l}(t)}{C_{k}^{-}(t)-\xi_{l}(t)}+\frac{\mu_{k}^{+}b_{l}(t)}{C_{k}^{+}(t)-\xi_{l}(t)}\right)
- \frac{1}{2}\left(\sum_{j=1}^{n}\frac{\alpha_{j}b_{l}(t)}{A_{j}(t)-\xi_{l}(t)}\right)\label{briggs}\end{equation}\normalsize
for all $l\in\{1,\ldots,N\}$. Where $C_{k}^{-}(t)$ and
$C_{k}^{+}(t)$ are the two preimages of $\xi_{k}(0)$ under $f_{t}$
satisfying $C_{k}^{-}(t)<C_{k}^{+}(t)$;
\[\mu_{k}^{\pm}=\deg_{H_{t},Q}(f_{t}(C_{k}^{\pm}(t)));\]
$A_{j}(t)=f_{t}^{-1}(a_{j})$; and
\[\Sigma_{0}=\left(\sum_{k=1}^{N}\beta_{k}\xi_{k}(0)\right) + \left(\sum_{j=1}^{n}\alpha_{j}a_{j}\right).\]
\end{theorem}
\begin{proof}
By Theorem \ref{QDSC} and Lemma \ref{degbound}, we can write
\[Q(w)=\left(\prod_{k=1}^{N}(w-\xi_{k}(0))^{\beta_{k}}\right)\left(\prod_{j=1}^{n}(w-a_{j})^{\alpha_{j}}   \right).\]
Then either by modifying the proof of Corollary \ref{QDSCslit} or iterating $N$ slit functions and applying Corollary \ref{QDSCslit} $N$ times, it is not too difficult to see that if we define $Q_{t}(z)$ by (\ref{calzaghe}), then
\begin{equation}Q_{t}(z)= \left(\prod_{k=1}^{N}(z-\xi_{k}(t))^2(z-C_{k}^{-}(t))^{\mu_{k}^{-}}(z-C_{k}^{+}(t))^{\mu_{k}^{+}}\right)\left(\prod_{j=1}^{n}(z-A_{j}(t))^{\alpha_{j}}\right).\label{manfredo}\end{equation}
Then the proof of (\ref{rahman}) is exactly the same as the proof of (\ref{toney}) in Theorem \ref{QDdriving}. The proof of (\ref{briggs}) is more complicated. First let
\[P_{t}(z)=-2\left(\sum_{k=1}^{N}\frac{b_{k}(t)}{z-\xi_{k}(t)}\right).\]
Then (\ref{urango}) becomes
\[\dot{f}_{t}(z)=f_{t}'(z)P_{t}(z).\]
Now take the logarithmic derivative of $Q_{t}(z)$ with respect to $z$ and $t$ separately using the definition of $Q_{t}(z)$ given by
(\ref{calzaghe}) to get
\[\frac{Q_{t}'(z)}{Q_{t}(z)} = \frac{Q'(f_{t}(z))f_{t}'(z)}{Q(f_{t}(z))}+2\frac{f_{t}''(z)}{f_{t}'(z)},\]
and
\begin{eqnarray*}\frac{\dot{Q}_{t}(z)}{Q_{t}(z)} &=& \frac{Q'(f_{t}(z))\dot{f}_{t}(z)}{Q(f_{t}(z))}+2\frac{\dot{f}_{t}'(z)}{f_{t}'(z)}\\ &=&
\frac{Q'(f_{t}(z))f_{t}'(z)P_{t}(z)}{Q(f_{t}(z))}+2\frac{f_{t}''(z)P_{t}(z)+f_{t}'(z)P_{t}'(z)}{f_{t}'(z)},
 \end{eqnarray*}
where we substitute (\ref{urango}) in for $\dot{f}_{t}$ to get from
the first to the second line. Thus we get
\begin{equation}\frac{\dot{Q}_{t}(z)}{Q_{t}(z)}=
\frac{Q_{t}'(z)}{Q_{t}(z)}P_{t}(z) +
2P_{t}'(z).\label{spinks}\end{equation}
 So then by (\ref{manfredo}), we note that
\[\frac{Q'_{t}(z)}{Q_{t}(z)}=\left(\sum_{k=1}^{N}\frac{2}{z-\xi_{k}(t)}+\frac{\mu_{k}^{-}}{z-C_{k}^{-}(t)} + \frac{\mu_{k}^{+}}{z-C_{k}^{+}(t)}\right)+\left(\sum_{j=1}^{n}\frac{\alpha_{j}}{z-A_{j}(t)}\right),
\]
and
\[\frac{\dot{Q}_{t}(z)}{Q_{t}(z)}=-\left[\left(\sum_{k=1}^{N}\frac{2\dot{\xi}_{k}(t)}{z-\xi_{k}(t)}+\frac{\mu_{k}^{-}\dot{C}_{k}^{-}(t)}{z-C_{k}^{-}(t)} + \frac{\mu_{k}^{+}\dot{C}_{k}^{+}(t)}{z-C_{k}^{+}(t)}\right)+\left(\prod_{j=1}^{n}\frac{\alpha_{j}\dot{A}_{j}(t)}{z-A_{j}(t)}\right)\right].
\]

\pagebreak
\begin{figure}[hp]
 \begin{center}
\scalebox{0.5}{\includegraphics{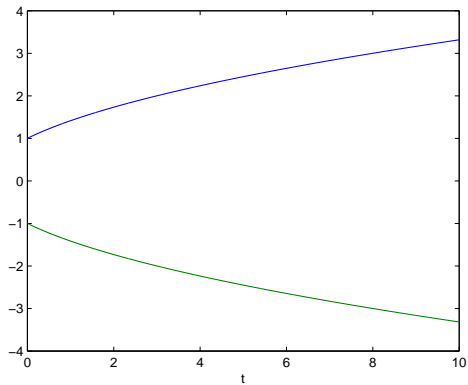}}
 \end{center}
\caption{A plot of the two driving functions $\xi_{1}$ (top)
and $\xi_{2}$ (bottom) on the $y$-axis against time on the
$x$-axis.} \label{fig4}
\end{figure}
Thus substituting this into (\ref{spinks}) and comparing the
coefficient of
\[\frac{1}{z-\xi_{l}(t)}\] (i.e. the residue at $z=\xi_{l}(t)$ of
both sides of (\ref{spinks})), we find that this is exactly
(\ref{briggs}).
\end{proof}

Similarly, we can prove a version of Theorem \ref{gammadot} and Corollary \ref{smoothxi} for multiple slits. This means that we can use the method
detailed in Section 4 with (\ref{briggs}) to calculate the driving
function for multiple $\theta_{k}$-trajectory arc slits. For example
Figure \ref{fig4} plots the graph of the two driving functions
$\xi_{1}$ and $\xi_{2}$ in the case when $\gamma_{1}$ and
$\gamma_{2}$ are 2 vertical slits starting from -1 and 1 (i.e.
orthogonal trajectories of $1dz^2$) and growing at the same speed.
Compare this with Figure 7 in \cite{KNK04}.
\subsection{Radial Loewner evolution}
The chordal Loewner differential equation was introduced because the
upper half-plane was an easier domain to work with for many
applications but the original setting of the Loewner differential
equation is in the unit disc $\mathbb{D}=\{z\in\mathbb{C}:|z|<1\}$: Suppose that $\gamma:[0,T)\mapsto
\overline{\mathbb{D}}$ is a simple curve such that $\gamma(0)\in\mathbb{T}=\{z:|z|=1\}$ and $\gamma(0,T)\subset\mathbb{D}\setminus\{0\}$. Then
$D_{t}=\mathbb{D}\setminus\gamma(0,t]$ is simply-connected and $0\in
D_{t}$ for all $t\in(0,T)$. Hence the Riemann mapping theorem implies that
there is unique conformal map $f_{t}$ mapping $\mathbb{D}$
conformally onto $D_{t}$ such that $f_{t}(0)=0$ and $f_{t}'(0)>0$. Then
Schwarz's lemma and the Carath\'{e}odory kernel theorem implies that $f_{t}'(0)$. is strictly decreasing and continuous so we
can reparameterize such that $f_{t}'(0)=e^{-t}$.  $f_{t}'(0)$ is sometimes called
the \emph{conformal radius} of $D_{t}$; hence in this case we are \emph{parameterizing by conformal radius}. Then the functions
$f_{t}$ satisfy the radial Loewner differential equation:
\[\dot{f}_{t}(z)=-zf_{t}'(z)\frac{z+e^{i\xi(t)}}{z-e^{i\xi(t)}}\]
See \cite{Low23} for more details.

\begin{theorem}
Suppose that $Q(w)dw^2$ is a K\"{u}hnau quadratic differential on
$\mathbb{D}$ such that $\deg_{Q}(0)=K\in\mathbb{Z}$ and
$e^{i\xi_{0}}\in\mathbb{T}=\{|z|=1\}$ satisfies
\[\deg_{\mathbb{D},Q}\left(e^{i\xi_{0}}\right)=N\in\{0,1,2,\ldots\}\]
then we have
\[Q(w)=w^K(w-e^{i\xi_{0}})^{N}\left(\prod_{j=1}^{n}(w-a_{j})^{\alpha_{j}}   \right),\]
where $a_{j}\in\mathbb{C}$ and $\alpha_{j}\in\mathbb{R}$.
Then if $\gamma:[0,T)\mapsto\overline{\mathbb{D}}$ is a simple curve such that $\gamma(0)=e^{i\xi_{0}}$ and $\gamma(0,T)\subset\mathbb{D}\setminus\{0\}$ is a $\phi$-trajectory arc
of $Q(w)dw^2$ in $\mathbb{D}$ that does not
meet $0$ and is parameterized
as above. Then we have
\begin{equation}e^{2i\xi(t)}= e^{-2t}\Pi_{0}C^{-}(t)^{-\mu^{-}}C^{+}(t)^{-\mu^{+}}\prod_{j=1}^{n}A_{j}(t)^{-\alpha_{j}},\label{mora}\end{equation}
and
\begin{equation}\dot{\xi}(t)=-\frac{1}{2i}\left( \mu^{-}\frac{C^{-}(t)+e^{i\xi(t)}}{C^{-}(t)-e^{i\xi(t)}}+\mu^{+}\frac{C^{+}(t)+e^{i\xi(t)}}{C^{+}(t)-e^{i\xi(t)}} + \sum_{j=1}^{n} \alpha_{j}\frac{A_{j}(t)+e^{i\xi(t)}}{A_{j}(t)-e^{i\xi(t)}}+2\right),\label{delahoya}\end{equation}
where, as usual, the functions $A_{j}(t)$ are defined by
\[A_{j}(t)= f_{t}^{-1}(a_{j}) \text{ for } j=1,\ldots,n,\]
and\[\mu^{\pm}=\deg_{D_{t},Q}(f_{t}(C^{\pm}(t))),\]
 $C^{+}(t)>C^{-}(t)$ are the two preimages of $e^{i\xi_{0}}$
under $f_{t}$; and also,
\[\Pi_{0}=e^{iN\xi_{0}}\prod_{j=1}^{n}a_{j}^{\alpha_{j}}.\]
\end{theorem}
\begin{proof}The formula for $Q(w)$ can be obtained from Theorem \ref{QDSC} by the transformation law. We then define $Q_{t}$ by (\ref{calzaghe}). Since the point 0 is fixed by
$f_{t}$, this implies that the $\deg_{Q_{t}}(0)=\deg_{Q}(0)=K$. Thus
we can apply Corollary \ref{QDSCslit} (again, using the transformation law) to get
\[Q_{t}(z)=z^K(z-e^{i\xi(t)})^2(z-C^{-}(t))^{\mu^{-}}(z-C^{+}(t))^{\mu^{+}}\left(\prod_{j=1}^{n}(z-A_{j}(t))^{\alpha_{j}} \right).\]
Then since by definition,
\[Q_{t}(z)=Q(f_{t}(z))f_{t}'(z)^{2}.\]
This immediately implies (\ref{mora}) by substituting $z=0$. Then we
get (\ref{delahoya}) in the same way as we get (\ref{wright}) from
(\ref{toney}) in the proof of Theorem \ref{QDdriving}.
\end{proof}
As in the case of multiple slits, a version of Theorem \ref{gammadot} and
Corollary \ref{smoothxi} holds for this case.
\subsection{Other versions of the Loewner differential equation}
There are several other versions of the Loewner differential
equation for simply-connected domains in the literature; the methods
in this paper should work in those cases as well and the proofs
should be similar to the proofs of Theorem \ref{QDdriving} etc. Also,
\cite{Kom43}, \cite{Kom50} generalizes the Loewner differential
equation to multiply-connected domains and again, some of the
methods should work in these cases possibly using methods in
\cite{Crow05} to extend Theorem \ref{QDSC} to multiply-connected
domains. Finally, even if we consider general 2-dimensional growth
processes given by the Loewner-Kufarev differential equation (see
Chapter 6 of \cite{Pom75}), some of the methods in this paper should
still be applicable.

\subsection*{Acknowledgement.} The author would like to express his
gratitude to his supervisor Dr. T. K. Carne for his constant
guidance and helpful discussion when writing this paper.

\end{document}